\documentclass[12pt]{article}
\usepackage{amsmath,amsthm,latexsym}
\usepackage{graphicx}
\newtheorem{theorem}{Theorem}
\newtheorem{remark}{Remark}

\newtheorem{lemma}{Lemma}
\newcommand{\h}{\hspace*{.24in}}
\def\geqslant {\ge}
\def\leqslant {\le}
\def\bq{\begin{equation}}
\def\eq{\end{equation}}
\def\bqq{\begin{eqnarray*}}
\def\eqq{\end{eqnarray*}}

\title{Determine the spacial term of a two-dimensional heat source}
\author{\normalsize Dang Duc Trong$^a$, ~Pham Ngoc Dinh Alain$^b$ and Phan Thanh Nam$^a$\\
\\\small\it $^a$Mathematics Department, HoChiMinh City National University, Viet Nam
\\\small\it $^b$Mathematics Department, Mapmo UMR 6628, BP 67-59, 45067 Orleans cedex, France
}
\date{{}}
\begin{document}
\maketitle
\begin{abstract} We consider the problem of determining a pair of functions $(u,f)$ satisfying the heat equation $u_t -\Delta u =\varphi(t)f (x,y)$, where $(x,y)\in \Omega=(0,1)\times (0,1)$ and the function $\varphi$ is given. The problem is ill-posed. Under a slight condition on $\varphi$, we show that the solution is determined uniquely from some boundary data and the initial temperature. Using the interpolation method and the truncated Fourier series, we construct a regularized solution of the source term $f$ from non-smooth data. The error estimate and numerical experiments are given.
\end{abstract}
\text{}\\
{\bf 1. Introduction}\\

Let $T>0$ and let $\Omega=(0,1)\times (0,1)$ be a heat conduction body. We consider the problem of determining a pair of functions $(u,f)$ satisfying the system
\bq
\left\{ \begin{gathered}
  u_t  - \Delta u = \varphi (t)f(x,y), \hfill \\
  u_x (0,y,t) = u_x (1,y,t) = u_y (x,0,t) = u_y (x,1,t) = 0, \hfill \\
  u(1,y,t) = 0 ,\hfill \\
u(x,y,0) =g(x,y),\hfill\\
 \end{gathered}  \right.\label{1}
\eq
for $(x,y)\in \Omega$, $t\in (0,T)$, where $g\in L^1(\Omega)$ and $\varphi\in L^1(0,T)$ are given.

This is a case of the problem of finding the source $F(\xi,t)$ satisfying the heat equation
\[
u_t  - \Delta u = F,
\]
where $\xi$ is the spacial variable. This inverse source problem is ill-posed. Indeed, a solution corresponding to the given data may not exist, and even if the solution exists (uniquely) then it may not depend continuously on the data. Because the problem is severely ill-posed and difficult, many presumptions on the form of the heat source are required. Roughly speaking, the function $F$ can be approximated by
\[
F(\xi ,t) \approx  \sum\limits_{n = 0}^N {\varphi _n (t)f_n (\xi )} .
\]
For simplicity, one reduces this approximation to its first term
$$F(\xi,t)=\varphi(t)f(\xi),$$
where one of the two functions $\varphi$ and $f$ is given. Many authors considered the uniqueness and stability conditions of the determination of the heat source  under this separate form \cite{CP90,CP91, Y93,Y94, YZ01,STY02,STY03,CY04}. In spire of the uniqueness and stability results, the regularization problem for unstable cases is still difficult. For a long time, it has been investigated for the heat souce which is time-depending only \cite{WZ99,SZN05,FL06} or space-depending only \cite{C68,WZ06,FL06}. Recently, the regularization problem for the heat source $F(\xi,t)=\varphi(t)f(\xi)$, where $\varphi$ is a given function, was regarded for one-dimensional case \cite{TLD05} and two-dimensional case \cite{TQD06}. However, these authors needed in addition an essential datum, that is the the final condition $u(\xi,T)$. Although this condition is unnatural, it gives an explicit representation of the solution as the inverse Fourier transform of a known term, and hence ones could use the truncated integral method to construct a regularized solution.

In the present  paper, we consider a similar problem to \cite{TLD05,TQD06}, but the final condition is removed completely. Moreover, the overspecified condition, i.e. $u(1,y,t)$, is slighter than this one in \cite{TQD06} and is almost optimal to still hold the uniqueness of the solution (see Remark 3). To our knowledge, no explicit form of the solution of system (\ref{1}) is available, and hence it is not easy to solve the problem although one has exact data. Of course, the problem with approximate data is even more difficult because of the ill-posedness.

Under a slight condition on $\varphi$, we shall  use the variational method and some properties of analytic functions to show the uniqueness of the solution. In particular, this result makes a regularization theorem of \cite{TQD06} trivial (see Remark 2). In spite of the uniqueness result, the problem is still ill-posed. We mention that the existence problem of a solution is not considered here. Instead, we shall assume that there is a (unique) exact solution corresponding to the exact data, and our aim is of constructing a regularized solution from approximate data. Using the interpolation method, we shall seek the  coefficients of the Fourier series expansion of the source term $f$ and then construct a regularized solution by the truncated Fourier series. The error estimate between the regularized solution and the exact solution is of order $(\ln(\varepsilon^{-1}))^{-1}$, where $\varepsilon$ is the error between the given data and the exact data. We also note that we shall concentrate only on finding the source term $f$ because we shall get a classical heat problem as soon as we know this function.

The remainder of the paper is divided into four sections.  We shall introduce some notations and state main results in Section 2. After that, we shall prove the uniqueness result in Section 3 and the regularization result in Section 4. In Section 5, we shall show how our method can be numerically implemented and give two examples to illuminate its effect.
\\\\
{\bf 2. Notations and main results}\\

By variational method, we have the following formula  to reconstruct the solution of the system $(1)$.
\begin{lemma}\label{Lem1} Assume that $(u,f)\in (C^1([0,T];L^1(\Omega))\cap L^2(0,T;H^2(\Omega)),L^2(\Omega))$
is a solution of the system $(1)$. Then for all $(\alpha,n)\in C\times Z$ we get
\bqq
&~& e^{-(\alpha ^2-n^2\pi^2)T }  \int\limits_{\Omega}{u(x,y,T)\cosh (\alpha x)\cos(n\pi y)dxdy}  \hfill \\
&~&-\int\limits_{\Omega} {g(x,y)\cosh (\alpha x)\cos(n\pi y)dxdy} \hfill \\
 &=& \int\limits_0^T {e^{-(\alpha ^2-n^2\pi^2)t} \varphi (t)dt} .\int\limits_{\Omega} {f(x,y)\cosh (\alpha x)\cos(n\pi y)dxdy}.
\eqq
\end{lemma}
From Lemma 1, we introduce some useful notations. For $w\in L^1(\Omega)$, $\varphi \in L^1(0,T)$ and $\alpha,\beta \in C$, put
\bqq
G(w)(\alpha,\beta ) &=& \int\limits_\Omega {w(x,y)\cosh(\alpha x)\cos(\beta y)dxdy},\\
D(\varphi)(\alpha,\beta) &=& \int\limits_0^T {e^{-(\alpha^2-\beta^2)t} \varphi(t)dt} ,\hfill\\
H(\varphi ,w)(\alpha,\beta )& = &\left\{ \begin{gathered}
  -\frac{{G(w)(\alpha,\beta )}} {{D(\varphi )(\alpha,\beta )}}, \h \text{ if }  D(\varphi )(\alpha,\beta )\ne 0, \hfill \\
  0,\h\h\h\h\h\text{if }  D(\varphi )(\alpha,\beta )=0.\hfill \\
\end{gathered}  \right.
\eqq
Note that for $w\in L^2(\Omega)$ and integers $m,n$,
\bq
G(w)(im\pi ,n\pi ) = \int\limits_\Omega  {w(x,y)\cos (m\pi x)\cos (n\pi y)dxdy}.\label{GF}
\eq
Because $\left\{ {\sqrt {\kappa (m,n)} \cos (m\pi x)\cos (n\pi y)} \right\}_{m\ge 0,n\ge 0}$ is an orthonormal basis on $L^2(\Omega)$, where
\[
\kappa (m,n) = \left\{ \begin{gathered}
  1, \text{ if } m=n=0, \hfill \\
  2, \text{ if } m>n=0 \text{ or } n>m=0, \hfill \\
  4, \text{ if } m>0 \text{ and } n>0,  \hfill \\
\end{gathered}  \right.
\]
we have the following representation
\bq
w(x,y) = \sum\limits_{m,n \geqslant 0} {\kappa (m,n)G(w)(im\pi,n\pi )\cos (m\pi x)\cos (n\pi y)}. \label{Fs}
\eq
This formula allows us to recover $f$ from $G(f)$. From Lemma 1, if $(\alpha^2-n^2\pi^2)>0$ is large and $|D(\varphi)(\alpha,n\pi)|$ is not so small then $G(f)(\alpha,n\pi)$ can be approximated by $H(\varphi,g)(\alpha,n\pi)$. To control $|D(\varphi)|$, we need the following condition (H) on $\varphi$.

(H) There exist $T_0\in(0,T]$, $\theta\ge0$ and $\Lambda>0$ and  such that either $\varphi(t)\ge \Lambda t^\theta$ for a.e $t\in (0,T_0)$, or $\varphi(t)\le -\Lambda t^\theta$ for a.e $t\in (0,T_0)$.

\begin{remark} The class of functions satisfying (H) is very broad. This condition holds with respect to $\theta=0$, for example, if $\varphi$ is continuous at $t=0$ and $\varphi(0)\ne 0$. To compare, we refer to the condition $\varphi\in C^1[0,T]$ and $\varphi(0)\ne 0$ in $\cite{ Y93,Y94}$.
\end{remark}

Under the condition $(H)$, we will obtain the uniqueness of the problem $(\ref{1})$.

\begin{theorem} Assume that $g\in L^1(\Omega)$ and $\varphi$ satifies $(H)$. Then the system $(\ref{1})$ has at most one solution $(u,f)$ in $\left( {C^1([0,T];L^1(\Omega))\cap L^2(0,T;H^2(\Omega)),L^2(\Omega)} \right).$
\end{theorem}

In spite of the uniqueness, the problem is still ill-posed and hence a regularization is necessary. Generally, the main ideas of the regularization are divided into three steps. For each integer $n$, we first approximate $G(f)(\alpha,n\pi)$ by $H(\varphi,g)(\alpha,n\pi)$ for some real numbers $\alpha$. In the next step we recover $G(f)(z,n\pi)$ when $z$ is in a ball of the complex  plane. Finally, we use a truncated series from the formula $(\ref{Fs})$ to construct the regularized solution.

For each integer $n$, we shall use the Lagrange interpolation polynomial to handle the key point of recovering $G(f)(.,n\pi)$. Recall that if $A=\{x_1,...,x_p\}$ be  a set of $p$ mutually distinct complex numbers and $w$ be a complex function then the Lagrange interpolation polynomial $L[A;w]$ is
\[
L[A;w](z) = \sum\limits_{j = 1}^p {\left( {\prod\limits_{k \ne j} {\frac{{z - x_k }}
{{x_j  - x_k }}} } \right)w(x_j )} .
\]
Now we are ready to state the regularization result.

\begin{theorem} Let $g_0\in L^1(\Omega)$ and let $\varphi_0\in L^1(0,T)$ satisfy $(H)$. Assume that $(u_0,f_0)\in (C^1([0,T];L^1(\Omega))\cap L^2(0,T;H^2(\Omega)),L^2(\Omega))$ is the exact solution of the system $(1)$ corresponding to the exact data $g_0$ and $\varphi_0$. Let $\varepsilon>0$, $\varphi_{\varepsilon}\in L^1(0,T)$ and $g_{\varepsilon}\in L^1(0,1))$ such that
\[
\left\| {\varphi_{\varepsilon}  - \varphi _0 } \right\|_{L^1 (0,T)}  \leqslant \varepsilon ,\left\| {g_{\varepsilon} - g_0 } \right\|_{L^1 (0,T)}  \leqslant \varepsilon .
\]
The regularized solution $f_\varepsilon$ is constructed from $\varphi_\varepsilon$ and $g_\varepsilon$ as follows
\bqq
r_\varepsilon&\in& Z\cap \left[{\frac{\ln (\varepsilon ^{ - 1} )}{50},\frac{\ln (\varepsilon ^{ - 1})}{50}+ 1} \right),\hfill\\
B(r_\varepsilon)&=& \{  \pm (4r_\varepsilon +j)|j = 1,2,...,20r_\varepsilon  \},\hfill\\
F_\varepsilon(m,n)&=&L\left[ {B(r_\varepsilon)  ;H(\varphi _\varepsilon  ,g_\varepsilon  )(.,n\pi )} \right](im\pi ),\hfill\\
f_\varepsilon  (x,y) &=& \sum\limits_{0 \leqslant m,n \leqslant r_\varepsilon  } {\kappa (m,n)F_\varepsilon(m,n)\cos (m\pi x)\cos (n\pi y)}.
\eqq
Then $f_\varepsilon\in C^{\infty}(R^2)$ and $\mathop {\lim }\limits_{\varepsilon  \to 0^ +  } f_\varepsilon   = f_0$ in $L^2(\Omega)$. Moreover, if $f_0\in H^1(\Omega)$ then $\mathop {\lim }\limits_{\varepsilon  \to 0^ +  } f_\varepsilon   = f_0$ in $H^1(\Omega)$ and there exists a constant $\varepsilon_0>0$ depending only on the exact data such that
\[
\left\| {f_\varepsilon   - f_0 } \right\|_{L^2 (\Omega )}  \leqslant \frac{{50}}
{\pi\ln (\varepsilon ^{ - 1} ) }\left\| {f_0 } \right\|_{H^1 (\Omega )}
\]
for all $\varepsilon\in (0,\varepsilon_0)$.
\end{theorem}
\text{}\\
{\bf 3. Uniqueness}
\\\\
{\bf Proof of Lemma 1.}
\begin{proof} Getting the inner product in $L^2(\Omega)$ of the first equation of the system $(1)$ and $W(x,y)=\cosh(\alpha x)\cos(n\pi y)$, then using the integral by part we have
\[
\frac{d}
{{dt}}\int\limits_\Omega  {uWdxdy}  - (\alpha ^2  - n^2 \pi ^2)\int\limits_\Omega  {uWdxdy} = \varphi (t)\int\limits_\Omega  {fWdxdy} .
\]
Next, we multiply the latter equality with $e^{-(\alpha ^2  -n^2 \pi ^2)t}$ to get
\bq
\begin{gathered}
\frac{d}
{{dt}}\left({e^{-(\alpha ^2  - n^2 \pi ^2 )t}  \int\limits_\Omega  {uWdxdy} } \right) = e^{-(\alpha ^2 - n^2 \pi ^2 )t} \varphi (t)\int\limits_\Omega  {fWdxdy} \label{tam}.\hfill \\
\end{gathered}
\eq
Finally, integrating $(\ref{tam})$ with respect to $t$ from 0 to $T$ we obtain the desired result.
\end{proof}
Now we consider some properties of the functions $G(w)$ and $D(\varphi)$.
\begin{lemma}\label{Lem2} Let $w\in L^1 (\Omega)$ and $n\in Z$. Then $G(w)(.,n\pi)$ is an entire function and
$$\left| {G(w)(z,n\pi )} \right| \leqslant e^{|z|}  \left\| w \right\|_{L^1 (\Omega)},~\forall z\in C.$$
Moreover, if $w\in L^2(\Omega)$ and $w\not  \equiv 0$ then there exists an integer $n$ such that
\[
\mathop {\lim \sup }\limits_{r \to  + \infty } \frac{{\ln |G(w)(r,n\pi )|}}
{r} \geqslant  - 1.
\]
\end{lemma}
\begin{proof} For each integer $n$, put  $\Phi_n(z)=G(w)(z,n\pi)$. Note that $z\mapsto \Phi_n(iz)$ is the cos-Fourier transform of the function
\[
x \mapsto \int\limits_0^1 {w(x,y)\cos (n\pi y)dy} ,~x \in [0,1].
\]
and hence $\Phi_n(iz)$ as well as $\Phi_n(z)$ are entire functions. Moreover,
\bqq
\left| {\Phi_n(z)} \right| &\le& \int\limits_{\Omega}{\left| {w(x,y)\cosh (zx)\cos(n\pi y)} \right|dxdy}\\
&\le& \int\limits_{\Omega}{\left| {w(x,y)\cosh (z)} \right|dxdy}
\le  e^{|z| } \left\| w \right\|_{L^1 (\Omega)},~\forall z\in C.
\eqq
Now assume that $w\not \equiv 0$. Note that
\[
\frac{d}
{{dz}}\Phi _n (im\pi) = \int\limits_\Omega  {ixw(x,y)\sin(m\pi x)\cos (n\pi y)dxdy}.
\]
and $\{\sin(m\pi x)\cos (n\pi y)\}_{m\ge 1, n\ge 0}$ is a orthogonal basis on $L^2(\Omega)$. Therefore, $\Phi_n$ is not constant for some integer $n$, and hence $ M_{\Phi_n}(r)>1$ for $r>0$ large enough, where $M_{\Phi_n} (r) = \mathop {\max }\limits_{|z| = r} \left| {\Phi_n(z)} \right|$. We shall show that
$$
\mathop {\lim \sup }\limits_{r \to  + \infty } \frac{{\ln|\Phi_n(r)|}}
{r} \ge -1.
$$
Of course, it is sufficient to consider the case $\mathop {\lim \sup }\limits_{r \to  + \infty } \ln |\Phi_n (r)|<0$. For $r>0$ large enough, since $\ln |\Phi(r)|<0$ and $ 1<M_{\Phi_n}(r)\le e^{r}\left\| w \right\|_{L^1 (\Omega)}$, we have
\[
 \frac{{\ln |\Phi_n(r)|}}
{{r+\ln \left\| w \right\|_{L^1 (\Omega )}}}\ge \frac{{\ln |\Phi_n(r)|}}
{{\ln M_{\Phi_n} (r)}}.
\]
Moreover, according to Beurling theorem (see, e.g., \cite{Ya}, Section 6.1, page 40) we get
\[
\mathop {\lim \sup }\limits_{r \to +\infty } \frac{{\ln \left| {\Phi_n(r)} \right|}}
{{\ln M_{\Phi_n}  (r)}} \geqslant  - 1.
\]
This implies the desired result.
\end{proof}
\begin{lemma}\label{Lem3} Let $\varphi\in L^1(0,T)$ and $(\alpha,n)\in R\times Z$. Then
$$\left| {D(\varphi )(\alpha,n\pi)} \right| \le ||\varphi||_{L^1(0,T)}$$
when $(\alpha^2-n^2\pi^2)>0$. Moreover, if $\varphi$ satisfies $(H)$ then
\[
\mathop {\lim \inf }\limits_{(\alpha ^2  - n^2 \pi ^2 ) \to  + \infty } (\alpha ^2  - n^2 \pi ^2 )^{\theta  + 1} \left| {D(\varphi )(\alpha ,n\pi )} \right| >0,
\]
where $\theta$ is as in $(H)$ corresponding to $\varphi$.
\end{lemma}
\begin{proof} The first inequality is obvious. Now assume that $\varphi$ satisfies the condition $(H)$ corresponding to $T_0$, $\theta$ and $\Lambda$. We shall prove that
\[
\mathop {\lim \inf }\limits_{\lambda  \to  + \infty } \lambda ^{\theta + 1} \left| {\int\limits_0^T {e^{ - \lambda t}\varphi(t)dt} } \right| >0,
\]
which will imply the desired result by choosing $\lambda=\alpha^2-n^2\pi^2$. We have
\bqq
  \lambda ^{\theta  + 1} \left| {\int\limits_0^T {e^{ - \lambda t} \varphi (t)dt} } \right|&\geqslant& \lambda ^{\theta  + 1} \left| {\int\limits_0^{T_0 } {e^{ - \lambda t} \varphi (t)dt} } \right| - \lambda ^{\theta  + 1} \left| {\int\limits_{T_0 }^T {e^{ - \lambda t} \varphi (t)dt} } \right| \hfill \\
   &\geqslant& \Lambda \lambda ^{\theta  + 1} \int\limits_0^{T_0 } {e^{ - \lambda t} t^\theta  dt}  - \lambda ^{\theta  + 1} e^{ - \lambda T_0 } \left\| \varphi  \right\|_{L^1 (0,T)}
\eqq
for all $\lambda>0$. Since $\mathop {\lim }\limits_{\lambda  \to  + \infty } \left( {\lambda ^{\theta  + 1} e^{ - \lambda T_0 } } \right) = 0$, it is sufficient to show that
\bq
\mathop {\lim \inf }\limits_{\lambda  \to  + \infty } \Psi _\theta  (\lambda ) > 0,\label{D1}
\eq
where
\[
\Psi _\theta  (\lambda ) = \lambda ^{\theta  + 1} \int\limits_0^{T_0 } {e^{ - \lambda t} t^\theta  dt} ,\h\forall \lambda>0.
\]
Using the integral by part we get
\[
\Psi _{\theta  + 1} (\lambda) =  - (\lambda T_0)^{\theta  + 1} e^{ - \lambda T_0 }  + (\theta  + 1)\Psi _\theta  (\lambda ).
\]
Therefore, it is enough to prove $(\ref{D1})$ for all $\theta\in [0,1)$. Indeed, by direct calculus we obtain $\mathop {\lim }\limits_{\lambda  \to  + \infty } \Psi _1 (\lambda ) = 1,\mathop {\lim }\limits_{\lambda  \to  + \infty } \Psi _2 (\lambda ) = 2$, and by Holder inequality we get
\[
\left( {\Psi _\theta  (\lambda )} \right)^{\frac{1}
{{2 - \theta }}} \left( {\Psi _2 (\lambda )} \right)^{\frac{{1 - \theta }}
{{2 - \theta }}}  \geqslant \Psi _1 (\lambda ), \h \forall\theta\in [0,1).
\]
Thus $\mathop {\lim \inf }\limits_{\lambda  \to  + \infty } \Psi _\theta  (\lambda ) \geqslant 2^{\theta  - 1}$ for all $\theta\in [0,1)$, and the proof is completed.
\end{proof}
\text{}\\
{\bf Proof of Theorem 1}
\begin{proof} Let $(u_1,f_1)$ and $(u_2,f_2)$ be two solutions of the system $(1)$. Put $u=u_1-u_2$ and $f=f_1-f_2$. We have to show that $(u,f)=(0,0)$.

Assume that $f\ne 0$. It follows from Lemma \ref{Lem1} that, for all $(\alpha,n) \in R\times Z$,
\bq
e^{-(\alpha^2-n^2\pi^2)T}G(u(.,.,T))(\alpha,n\pi)=D(\varphi)(\alpha,n\pi).G(f)(\alpha,n\pi).\label{1.1}
\eq
Due to Lemma \ref{Lem2}, there exist an integer $n$ and a sequence of positive numbers $\{\alpha_m\}$, which depends on $n$, such that $\mathop {\lim }\limits_{m \to \infty } \alpha _m  =  + \infty$ and, for all $m\ge 1$,
$$\left| {G(f)(\alpha _m,n\pi )} \right| \geqslant e^{ - 2\alpha_m }.$$
Using Lemma \ref{Lem2} again  we get
\[
|G(u(.,.,T))(\alpha _m,n\pi )| \leqslant e^{\alpha _m } \left\| {u(.,.,T)} \right\|_{L^1 (\Omega)}.
\]
Moreover, according to Lemma \ref{Lem3}, there exist constants $\theta>0$ and $C_0>0$, which depend on $\varphi$, such that
\[
\left| {D(\varphi )(\alpha _m,n\pi )} \right| \geqslant \frac{{C_0 }}
{{(\alpha _m^2-n^2\pi^2)^{\theta + 1} }}
\]
for all $m$ large enough. Therefore, from (\ref{1.1}) we obtain that
\[
e^{-(\alpha _m^2-n^2\pi^2)T}e^{\alpha_m} \left\| {u(.,.,T)} \right\|_{L^1 (\Omega)}  \geqslant \frac{{C_0}}
{{(\alpha _m^2-n^2\pi^2)^{\theta +1}}}e^{ -2\alpha _m}
\]
for all $m$ large enough. This is a contradiction.
\\\h Thus $f= 0$. Hence, the equality $(\ref{tam})$ in the proof of Lemma 1 reduces to
\[
\frac{d}{dt}\left( {e^{-(\alpha ^2-n^2\pi^2)t}G(u(.,.,t))(\alpha,n\pi) } \right) = 0.
\]
It follows from $u(.,.,0)= 0$ that $u=0$.
\end{proof}
\begin{remark} In Theorem 2 of \cite{TQD06}, the authors imposed that $u(.,.,0)= 0$ and $\varphi<-C_0$ a.e. for some positive constant $C_0$. However, due to our uniqueness result, these requirements are so strict that the problem has only the trivial solution.
\end{remark}
\begin{remark} In comparison to \cite{TQD06}, here the datum $u(x,1,t)$ is omitted. Moreover, if $u(1,y,t)$ is omitted too then the uniqueness of the solution may not hold even if the final temperature is given. For example, the system
$$
\left\{ \begin{gathered}
  u_t  - \Delta u =(\pi cos(\pi t)+2\pi^2 \sin(\pi t))f(x,y), \hfill \\
  u_x (0,y,t) = u_x (1,y,t) = u_y (x,0,t) = u_y (x,1,t) = 0,\hfill\\
u(x,y,0) = 0,u(x,y,1) = 0, \hfill \\
\end{gathered}  \right.
$$
has the trivial solution $(u,f)=(0,0)$ and (at least) a non-trivial solution
\[
u(x,y,t) = \sin (\pi t)\cos (\pi x)\cos (\pi y),f(x,y) = \cos (\pi x)\cos (\pi y).
\]
Thus the overspecified condition $u(1,y,t)$ is essential.
\end{remark}
\text{}\\
\text{\bf 4. Regularization}\\\\
\h To construct the regularized solution, we first want to solve the problem of recovering the entire function $G(f)(.,n\pi)$ for each integer $n$. The key tool is the Lagrange interpolation polynomial.
\begin{lemma}\label{Lem4} Let $r\ge 1$ be an integer and $B_r =\{\pm (4r+j)|j=1,2,...,20r\}$. Let $w$ and $\widetilde w$ be two even complex function such that $w$ is an entire function and $|w(z)|\le Ae^{|z|}$ for all $z\in C$, where $A$ is independent on $z$. Then
\[
\mathop {\sup }\limits_{|z|\le \pi r}\left| {w(z) - L(B_r ;\widetilde w)(z)} \right| \leqslant Ae^{ - r}  + 20re^{25r}\mathop {\sup }\limits_{z\in B_r} |w(z) - \widetilde w(z)|.
\]
\end{lemma}
\begin{proof} Fix $z\in C, |z|\le \pi r$ and denote $z_j=4r+j$ for each $j=1,2,...,20r$. We shall use the triangle inequality
\bqq
|w(z) - L(B_r;\widetilde w )(z)|\le | w(z) - L(B_r ;w)(z)|+|L(B_r ;w-\widetilde w)(z)|.
\eqq
\h We first estimate $|w(z) - L(B_r ;w)(z)|$. Let $\gamma=\{z\in C, |z|=45r\}$. Using the residue theorem, we get the Hermite's remainder formula
\[
w(z) - L[B_r;w](z) = \frac{1}
{{2\pi i}}\int\limits_\gamma  {\frac{{w(\xi )}}
{{\xi  - z}} \cdot \prod\limits_{j = 1}^{20r} {\frac{{z^2 - z^2_j }}
{{\xi^2  - z^2_j }}} d\xi } .
\]
Therefore,
\bq
\left| {w(z) - L[B_r;w](z)} \right| \leqslant 45r\mathop {\sup }\limits_{\xi  \in \gamma } \left\{ {\frac{{|w(\xi )|}}
{{|\xi  - z|}} \cdot \prod\limits_{j = 1}^{20r}{\frac{{|z^2 - z^2_j |}}
{{|\xi^2  - z^2_j |}}} } \right\}.\label{L2}
\eq
For $\xi\in \gamma$ we have $|w(\xi )| \leqslant Ae^{45r} $, $|\xi  - z| \geqslant (45 - \pi )r$ and
\bq
\prod\limits_{j = 1}^{20r} {\frac{{|z^2  - z_j^2 |}}
{{|\xi ^2  - z_j^2 |}}}  \leqslant \prod\limits_{j = 1}^{20r} {\frac{{|z|^2  + z_j^2 }}
{{|\xi |^2  - z_j^2 }}}  \leqslant \prod\limits_{j = 1}^{20r} {\frac{{(\pi r)^2  + z_j^2 }}
{{(45r)^2  - z_j^2 }}} .\label{L3}
\eq
We shall show that
\bq
\prod\limits_{j = 1}^{20r} {\frac{{(\pi r)^2  + z_j^2 }}
{{(45r)^2  - z_j^2 }}}  \leqslant \frac{45-\pi}{45}e^{ - 46r},~\forall r=1,2,...\label{L4}
\eq
We can check (\ref{L4}) by direct computations for   $r=1,2,...,54$. Now we consider when $r\ge 55$. Since the real function
\[
\vartheta (x) = \ln \left( {\frac{{(\pi r)^2 + x}}
{{(45r)^2 - x}}} \right).
\]
is increasing and concave in $[0,(24r)^2]$, we can apply Jensen's inequality to get
\bqq
   \sum\limits_{j = 1}^{20r} {\vartheta (z_j^2 )}
    & =& \sum\limits_{k = 1}^5{\left( {\sum\limits_{j = 4(k - 1)r + 1}^{4kr} {\vartheta (z_j^2 )} } \right)}\leqslant 4r\sum\limits_{k = 1}^5 {\vartheta \left( {\frac{1}
{{4r}}\sum\limits_{j = 4(k - 1)r + 1}^{4k} {z_j^2 } } \right)}  \hfill \\
     &=& 4r\sum\limits_{k = 1}^5 {\vartheta \left( {\left( {16k^2  + 16k + \frac{{16}}
{3}} \right)r^2  + (4k + 2)r + \frac{1}
{6}} \right)}  \hfill \\
   &\leqslant& 4r\sum\limits_{k = 1}^5 {\vartheta \left( {\left( {16k^2  + 16k + \frac{{16}}
{3}} \right)r^2  + (4k + 2)\frac{{r^2 }}
{{55}} + \frac{{r^2 }}
{{6 \times 55^2 }}} \right)}  \hfill \\
   &=& 4r\sum\limits_{k = 1}^5 {\vartheta \left( {\left( {16k^2  +\frac{884}{55}k + \frac{32487}{6050}} \right)r^2 } \right)}  \hfill \\
  &=& 4r\sum\limits_{k = 1}^5 {\ln \left( {\frac{{\pi ^2  + \left( {16k^2  +\frac{884}{55}k + \frac{32487}{6050}} \right)}}
{{45^2  - \left( {16k^2  +\frac{884}{55}k + \frac{32487}{6050}} \right)}}} \right)}  \hfill \\
   &=& 4r \times (-11.51809713)< - 46r+\ln \left( {\frac{{45 - \pi }}
{{45}}} \right),~r\ge 55.
\eqq
Thus $(\ref{L4})$ holds for all $r=1,2,...$. From $(\ref{L3})$ and $(\ref{L4})$, we reduce $(\ref{L2})$ to
\bq
|w(z) - L(B_r ;w)(z)| \le Ae^{ - r},~\forall r=1,2,...\label{L7}
\eq
\h We shall now estimate $|L(B_r;w-\widetilde w)(z)|$. Since
\[
L(B_r ;w - \widetilde w)(z) = \sum\limits_{j = 1}^{20r} {\left( {\prod\limits_{k \ne j} {\frac{{z^2  - z_k^2 }}
{{z_j^2  - z_k^2 }}} } \right)} \left( {w(z_j ) - \widetilde w(z_j )} \right),
\]
we obtain
\bq
\left| {L(B_r ;w - \widetilde w)(z)} \right| \leqslant 20r\sigma \mathop {\sup }\limits_{1 \leqslant j \leqslant 20r} \prod\limits_{k \ne j} {\frac{{|z|^2  + z_k^2 }}
{{|z_j^2  - z_k^2 |}}} ,\label{L5}
\eq
where $\sigma=\mathop {\sup }\limits_{1 \leqslant j \leqslant 20r} |w(z_j ) - \widetilde w(z_j )|$. We have
\bqq
  \prod\limits_{k \ne j} {\frac{{|z|^2  + z_k^2 }}
{{\left| {z_j^2  - z_k^2 } \right|}}}
&=& \prod\limits_{k \ne j} {\frac{{|z|^2  + z_k^2 }}
{{(z_j  + z_k )z_k }}}  \times \frac{{z_k }}
{{\left| {z_j  - z_k } \right|}} \leqslant \prod\limits_{k \ne j} {\frac{{z_k }}
{{\left| {z_j  - z_k } \right|}}}  \hfill \\
 &\leqslant& \frac{{\prod\limits_{k = 2}^{20r} {z_k } }}
{{\prod\limits_{k \ne j} {\left| {z_j  - z_k } \right|} }} = \frac{{(4r + 2).(4r+3)... (24r)}}
{{(j - 1)!(20r - j)!}} \hfill \\
 &\leqslant& \frac{{(4r + 2).(4r+3)... (24r)}}
{{(10r - 1)!(10r)!}} = :J(r).
\eqq
Note that $J(1)<e^{25}$ and
\bqq
  \frac{{J(r + 1)}}
{{J(r)}}
&=& \frac{{(24r+1).(24r +2)...(24r + 24)}}
{{(4r + 2)...(4r +5) \times (10r) \times \left[ {(10r + 1)...(10r + 9)} \right]^2  \times (10r + 10)}}+ \hfill \\
   &<& \frac{{24^{24} }}
{{4^4  \times 10^{20} }} < e^{25} .
\eqq
Here, the latter inequality can checked by direct expansion
\[
\begin{gathered}
  ~~24^{24}  \times (4r + 2)...(4r + 5) \times (10r) \times \left[ {(10r + 1)...(10r + 9)} \right]^2  \times (10r + 10) \hfill \\
   - 4^4  \times 10^{20}  \times (24r + 1)(24r + 2)...(24r + 24) \hfill \\
\end{gathered}
\]
to get a polynomial of degree 23 in term of $r$, which is obviously positive for $r\ge 1$.
Thus $J(r)<e^{25r}$ for each integer $r\ge 1$. Therefore, we can reduce $(\ref{L5})$ to
\bq
\left| {L(B_r ;w - \widetilde w)(z)} \right| \leqslant 20re^{25r}\sigma.\label{L6}
\eq
From (\ref{L7}) and (\ref{L6}) we have the desired result.
\end{proof}
In our application of Lemma 4, we shall choose $w=G(f_0)(.,n\pi)$, an unknown function, and $\widetilde w=H(\varphi_\varepsilon,g_\varepsilon)(.,n\pi)$, which is known from data. We have the following error estimate  between these functions.
\begin{lemma}\label{Lem5} Let $u_0,f_0,\varphi_0,g_0,\varphi_\varepsilon,g_\varepsilon,r_\varepsilon,B(r_\varepsilon)$ be as in Theorem 2 and $\theta$ be as in the condition $(H)$ corresponding to $\varphi_0$. Then there exists $\varepsilon_1>0$ depending only on the exact data such that
\[
\left| {G(f_0 )(\alpha ,n\pi ) - H(\varphi _\varepsilon  ,g_\varepsilon  )(\alpha ,n\pi )} \right| \leqslant \left( {\ln (\varepsilon ^{ - 1} )} \right)^{4\theta  + 5} e^{|\alpha |} \varepsilon .
\]
provided $\varepsilon\in (0,\varepsilon_1)$, $0\le n\le r_\varepsilon$ and $\alpha\in B(r_\varepsilon)$.
\end{lemma}
\begin{proof} Note that if $\alpha\in B(r_\varepsilon)$ then $4r_\varepsilon\le |\alpha|\le 24r_\varepsilon$. Hence, for $\varepsilon>0$ small enough one has
$$\frac{\ln (\varepsilon ^{ - 1} )}{T}\leqslant \alpha ^2  - n^2 \pi ^2  \leqslant \left( {\ln (\varepsilon ^{ - 1} )} \right)^2.$$
Thus, according to Lemma \ref{Lem3}, there exists $C(\varphi_0)>0$ depending only on $\varphi_0$ such that
\[
\left| {D(\varphi _0 )(\alpha ,n\pi )} \right| \geqslant \frac{{C(\varphi_0 )}}
{{(\alpha ^2  - n^2 \pi ^2 )^{\theta  + 1} }} \geqslant \left( {\ln (\varepsilon ^{ - 1} )} \right)^{ - 2(\theta  + 1)},
\]
and consequently,
\bqq
\left| {D(\varphi _\varepsilon  )(\alpha ,n\pi )} \right| &\ge& \left| {D(\varphi _0 )(\alpha ,n\pi )} \right| - \left| {D(\varphi _0 )(\alpha ,n\pi ) - D(\varphi _\varepsilon  )(\alpha ,n\pi )} \right|  \\
&\ge& \left( {\ln (\varepsilon ^{ - 1} )} \right)^{ - 2(\theta  + 1)}  - \varepsilon  \geqslant \frac{1}
{2}\left( {\ln (\varepsilon ^{ - 1} )} \right)^{ - 2(\theta  + 1)}.
\eqq
It follows from Lemma \ref{Lem1} and Lemma \ref{Lem2}, for $\varepsilon>0$ small enough, that
\[
\begin{gathered}
 ~~ \left| {G(f_0 )(\alpha ,n\pi ) - H(\varphi _0 ,g_0 )(\alpha ,n\pi )} \right| = e^{ - (\alpha ^2  - n^2 \pi ^2 )T} .\left| {\frac{{G(u_0 (.,.,T))(\alpha ,n\pi )}}
{{D(\varphi _0 )(\alpha ,n\pi )}}} \right| \hfill \\
   \leqslant e^{ - (\alpha ^2  - n^2 \pi ^2 )T} .\frac{{e^{|\alpha |} \left\| {u_0 (.,.,T)} \right\|_{L^1 (\Omega )} }}
{{\left( {\ln (\varepsilon ^{ - 1} )} \right)^{ - 2(\theta  + 1)} }} \leqslant \frac{1}{2}\varepsilon.e^{|\alpha |} \left( {\ln (\varepsilon ^{ - 1} )} \right)^{4\theta  + 5}. \hfill \\
\end{gathered}
\]
Moreover, for $\varepsilon>0$ small enough,
\[
\begin{gathered}
 ~~ \left| {H(\varphi _0 ,g_0 )(\alpha ,n\pi ) - H(\varphi _\varepsilon  ,g_\varepsilon  )(\alpha ,n\pi )} \right| \hfill \\
   = \left| {\frac{{G(g_0 (.,.,T))(\alpha ,n\pi )}}
{{D(\varphi _0 )(\alpha ,n\pi )}} - \frac{{G(g_\varepsilon  (.,.,T))(\alpha ,n\pi )}}
{{D(\varphi _\varepsilon  )(\alpha ,n\pi )}}} \right| \hfill \\
   \leqslant \frac{{\left| {G(g_0 (.,.,T))(\alpha ,n\pi )} \right|.\left| {D(\varphi _\varepsilon  )(\alpha ,n\pi ) - D(\varphi _0 )(\alpha ,n\pi )} \right|}}
{{\left| {D(\varphi _0 )(\alpha ,n\pi )D(\varphi _\varepsilon  )(\alpha ,n\pi )} \right|}} \hfill \\
 ~~  + \frac{{\left| {D(\varphi _0 )(\alpha ,n\pi )} \right|.\left| {G(g_\varepsilon  (.,.,T))(\alpha ,n\pi ) - G(g_0 (.,.,T))(\alpha ,n\pi )} \right|}}
{{\left| {D(\varphi _0 )(\alpha ,n\pi )D(\varphi _\varepsilon  )(\alpha ,n\pi )} \right|}} \hfill \\
   \leqslant \frac{{e^{|\alpha |} \left\| {g_0 } \right\|_{L^1 (\Omega )} .\varepsilon }}
{{\left( {\ln (\varepsilon ^{ - 1} )} \right)^{ - 2(\theta  + 1)} .2\left( {\ln (\varepsilon ^{ - 1} )} \right)^{ - 2(\theta  + 1)} }} + \frac{{\left\| {\varphi _0 } \right\|_{L^1 (0,T)} .e^{|\alpha |} \varepsilon }}
{{\left( {\ln (\varepsilon ^{ - 1} )} \right)^{ - 2(\theta  + 1)} .2\left( {\ln (\varepsilon ^{ - 1} )} \right)^{ - 2(\theta  + 1)} }} \hfill \\
   \leqslant \frac{1}{2}\left( {\ln (\varepsilon ^{ - 1} )} \right)^{4\theta  + 5} e^{|\alpha |} \varepsilon.  \hfill \\
\end{gathered}
\]
The desired result follows the two latter inequalities and the triangle inequality.
\end{proof}
For $w\in L^2(\Omega)$ and an integer $M\ge 1$, put
\[
\Gamma _M (w)(x,y) = \sum\limits_{0 \leqslant m,n \leqslant M} {\kappa (m,n)G(im\pi ,n\pi )\cos (m\pi x)\cos (n\pi y)}.
\]
The formula $(\ref{Fs})$ means $\mathop {\lim }\limits_{M \to  + \infty } \Gamma _M (w) = w$ in $L^2(\Omega)$. However, to prove Theorem 2 we shall need a sharper estimate for the remainder of the Fourier series.
\begin{lemma}\label{Lem6} If $w\in H^1(\Omega)$ then $\mathop {\lim }\limits_{M \to  + \infty } \Gamma _M (w) = w$ in $H^1(\Omega)$ and
\[
\left\| {\Gamma _M (w) - w} \right\|_{L^2 (\Omega )}  \leqslant \frac{1}
{{\pi (M+1)}}\left\| w \right\|_{H^1 (\Omega )} .
\]
\end{lemma}
\begin{proof} From $\ref{GF}$, using the integral by part we get
\bqq
  \pi m G(w)(im\pi ,n\pi ) &=&  - \int\limits_\Omega  {w_x (x,y)\sin (m\pi x)\cos (n\pi y)dxdy},  \hfill \\
  \pi n G(w)(im\pi ,n\pi ) &=&  - \int\limits_\Omega  {w_y (x,y)\cos (m\pi x)\sin (n\pi y)dxdy},
\eqq
and  hence
\bqq
&~&\pi ^2 (m^2  + n^2 )\kappa (m,n)\left| {G(w)(im\pi ,n\pi )} \right|^2 \hfill\\
&=&\kappa (m,n) \left( {\int\limits_\Omega  {w_x (x,y)\sin (m\pi x)\cos (n\pi y)dxdy} } \right)^2 \\
&~&+\kappa (m,n)\left( {\int\limits_\Omega  {w_y (x,y)\cos (m\pi x)\sin (n\pi y)dxdy} } \right)^2.
 \eqq
Since $\left\{ {\sin (m\pi x)\cos (n\pi y)} \right\}_{m\ge 1,n\ge 0}$ and $\left\{ {\cos (m\pi x)\sin (n\pi y)} \right\}_{m\ge 0,n\ge 1}$ are orthogonal bases on $L^2(\Omega)$, it follows from Parseval equality that
\[
\begin{gathered}
  \sum\limits_{m,n \geqslant 0} {\kappa (m,n)\left( {\int\limits_\Omega  {w_x (x,y)\sin (m\pi x)\cos (n\pi y)} } \right)^2 }  = \left\| {w_x } \right\|_{L^2 (\Omega )}^2 , \hfill \\
  \sum\limits_{m,n \geqslant 0} {\kappa (m,n)\left( {\int\limits_\Omega  {w_y (x,y)\cos (m\pi x)\sin (n\pi y)} } \right)^2 }  = \left\| {w_y } \right\|_{L^2 (\Omega )}^2 . \hfill \\
\end{gathered}
\]
Thus
\bqq
&~&\sum\limits_{m,n \geqslant 0} { \pi ^2(m^2  + n^2 )\kappa (m,n)\left| {G(w)(im\pi ,n\pi )} \right|^2 }\hfill\\
&=& \left\| {w_x } \right\|_{L^2 (\Omega )}^2  + \left\| {w_y } \right\|_{L^2 (\Omega )}^2  \le \left\| w \right\|_{H^1 (\Omega )}^2 .
\eqq
Using Parseval equality again, from the latter inequality we have
\bqq
  &~&\left\| {\Gamma _M (w)- w} \right\|_{L^2 (\Omega )}^2\hfill\\
&=&\sum\limits_{\max \{ m,n\}  > M} {\kappa (m,n)\left| {G(w)(im\pi ,n\pi )} \right|^2 }  \hfill \\
  &<&\frac{1}{{\pi ^2 (M+1)^2 }}\sum\limits_{\max \{ m,n\}  > M} {\pi ^2 (m^2  + n^2 )\kappa (m,n)\left| {G(w)(im\pi ,n\pi )} \right|^2 } \hfill \\
&\leqslant&\frac{1}{{\pi ^2 (M+1)^2 }}\left\| w \right\|_{H^1 (\Omega )}^2  ,
\eqq
and
\bqq
 &~&\left\| {\Gamma _M (w) - w} \right\|_{H^1 (\Omega )}^2\hfill\\
&=&\sum\limits_{\max\{m,n\}>M} {\left( {1 + \pi ^2 (m^2  + n^2 )} \right)\kappa (m,n)\left| {G(w)(im\pi ,n\pi )} \right|^2 } \to 0
\eqq
as $M\to +\infty$.
\end{proof}
\text{}\\
{\bf Proof of Theorem 2.}
\begin{proof} We shall first get the error between $\Gamma_{r_\varepsilon}(f_0)$ and $f_\varepsilon$, and then use the approximation between $\Gamma_{r_\varepsilon}( f_0)$ and $f_0$ to get the results. Note that
\[
G(f_\varepsilon  )(im\pi ,n\pi ) = \left\{ \begin{gathered}
  F_\varepsilon  (m,n) = L\left[ {B_\varepsilon  ;H(\varphi _\varepsilon  ,g_\varepsilon  )(.,n\pi )} \right](im\pi ),~\text{if}~ 0\le m,n\le r_\varepsilon, \hfill \\
  0,~\text{ otherwise.} \hfill \\
\end{gathered}  \right.
\]

Let us consider $0\le m,n\le r_\varepsilon$. By Lemma \ref{Lem2}, $G(f_0 )(. ,n\pi )$ is an entire function and
$\left| {G(f_0 )(z,n\pi )} \right| \leqslant \left\| {f_0 } \right\|_{L^1 (\Omega )} e^{|z|} ,~z \in C.$ Moreover, for $\varepsilon>0$ small enough, Lemma \ref{Lem5} gives
\[
\mathop {\sup }\limits_{\alpha  \in B(r_\varepsilon  )} \left| {G(f_0 )(\alpha ,n\pi ) - H(\varphi _\varepsilon  ,g_\varepsilon  )(\alpha ,n\pi )} \right| \leqslant (\ln (\varepsilon ^{ - 1} ))^{4\theta  + 5} e^{24r_\varepsilon  } \varepsilon .
\]
Therefore, it follows from Lemma \ref{Lem4} that
\[
\begin{gathered}
~~\left| {G(f_0 )(im\pi ,n\pi ) - G(f_\varepsilon  )(im\pi ,n\pi )} \right|\hfill\\
= \left| {G(f_0 )(im\pi ,n\pi ) - L\left[ {B_\varepsilon  ;H(\varphi _\varepsilon  ,g_\varepsilon  )(.,n\pi )} \right](im\pi )} \right| \hfill \\
   \leqslant \left\| {f_0 } \right\|_{L^1 (\Omega )} e^{ - r_\varepsilon  }  + 20r_\varepsilon  e^{25r_\varepsilon  }\left( {\ln (\varepsilon ^{ - 1} )} \right)^{4\theta  + 5} e^{24r_\varepsilon} \varepsilon, ~0\le m,n\le r_\varepsilon,\hfill \\
\end{gathered}
\]
for $\varepsilon>0$ small enough. Note that $e^{r_\varepsilon}\le e\varepsilon^{-1/50}$, the above inequality reduces to
$$
\left| {G(f_0 )(im\pi ,n\pi ) - G(f_\varepsilon  )(im\pi ,n\pi )} \right| \le \varepsilon^{1/51},~0\le m,n\le r_\varepsilon,
$$
for $\varepsilon>0$ small enough. Therefore,
\[
\begin{gathered}
 ~~ \left\| {\Gamma _{r_\varepsilon  } (f_0 ) - f_\varepsilon  } \right\|_{H^1 (\Omega )}^2  \hfill \\
   = \sum\limits_{0 \leqslant m,n \leqslant r_\varepsilon  } {\left( {1 + \pi ^2 (m^2  + n^2 )} \right)\kappa (m,n)\left| {G(f_0 )(im\pi ,n\pi ) - G(f_\varepsilon  )(im\pi ,n\pi )} \right|^2 }  \hfill \\
   \leqslant (1 + r_\varepsilon  )^2 \left( {1 + 2\pi ^2 r_\varepsilon ^2 } \right)4\varepsilon ^{2/51}  \leqslant \varepsilon ^{1/26}  \hfill \\
\end{gathered}
\]
for $\varepsilon>0$ small enough.

Thus $\mathop {\lim }\limits_{\varepsilon  \to 0^ +  } \left( {\Gamma _{r_\varepsilon  } (f_0 ) - f_\varepsilon  } \right) = 0$ in $H^1(\Omega)$. Hence it follows from $\mathop {\lim }\limits_{\varepsilon  \to 0^ +  } \Gamma _{r_\varepsilon  } (f_0 ) = f_0$ in $L^2(\Omega)$ that $\mathop {\lim }\limits_{\varepsilon  \to 0^ +  } f_\varepsilon   = f_0$ in $L^2(\Omega)$.

Now assume in addition that $f_0\in H^1(\Omega)$. Then Lemma \ref{Lem6} leads to $\mathop {\lim }\limits_{\varepsilon  \to 0^ +  } \Gamma _{r_\varepsilon  } (f_0 ) = f_0$ in $H^1(\Omega)$ , and hence $\mathop {\lim }\limits_{\varepsilon  \to 0^ +  } f_\varepsilon   = f_0$ in $H^1(\Omega)$. Moreover, using Lemma 6 again we get
\bqq
  \left\| {f_\varepsilon   - f_0 } \right\|_{L^2 (\Omega )}
&\leqslant& \left\| {\Gamma _{r_\varepsilon  } (f_0 ) - f_\varepsilon  } \right\|_{L^2 (\Omega )}  + \left\| {\Gamma _{r_\varepsilon  } (f_0 ) - f_0 } \right\|_{L^2 (\Omega )}  \hfill \\
  &\leqslant& \varepsilon ^{1/52}  + \frac{1}
{{\pi (r_\varepsilon   + 1)}}\left\| {f_0 } \right\|_{H^1 (\Omega )}  \hfill\\
&\leqslant& \frac{50}{{\pi \ln(\varepsilon^{-1})  }}\left\| {f_0 } \right\|_{H^1 (\Omega )}
\eqq
for $\varepsilon>0$ small enough.
\end{proof}
\text{}\\
{\bf 5. Numerical experiments}\\\\
\h In this section, for simplicity we shall choose $T=1$ and $\varepsilon=k^{-1}$, where $\varepsilon$ is the error of the data and $k$ is an integer. The regularized scheme in Theorem 2 can be rewritten as the following procedure, where $k$, $\varphi_k$, $g_k$ are given data and $f_k$ is the regularized solution.\\
\[
\begin{gathered}
   r: = {\text{ceil}}\left( {\frac{{\ln (k)}}
   {{50}}} \right);f_k: = 0;\hfill\\
  for{\text{ }}n{\text{ }}from{\text{ }}1{\text{ }}to{\text{ }}20r{\text{ }}do \hfill \\
  \h for{\text{ }}j{\text{ }}from{\text{ }}1{\text{ }}to{\text{ }}20r{\text{ }}do \hfill \\
  \h\h z_j : = 4r + j;H_j : =  - \frac{{\int\limits_0^1 {\int\limits_0^1 {g_k \cosh (z_j x)\cos (n\pi y)dxdy} } }}
{{\int\limits_0^T {e^{ - (z_j^2  - n^2 \pi ^2 )t} \varphi _k (t)dt} }}; \hfill \\
  \h end{\text{ }}do; \hfill \\
  \h for{\text{ }}m{\text{ }}from{\text{ }}0{\text{ }}to{\text{ }}r{\text{ }}do \hfill \\
  \h\h {\text{coef}}: = 0; \hfill \\
  \h\h for{\text{ }}j{\text{ }}from{\text{ }}1{\text{ }}to{\text{ }}20r{\text{ }}do \hfill \\
  \h\h\h w: = 1; \hfill \\
  \h\h\h for{\text{ }}\ell{\text{ }}from{\text{ }}1{\text{ }}to{\text{ }}20r{\text{ }}do \hfill \\
  \h\h\h\h if{\text{ }}(\ell\ne j){\text{ }}then{\text{ }} w: = w*\frac{{ - (m\pi )^2  - z_\ell^2 }}
{{z_j^2  - z_\ell^2 }}{\text{ }}end{\text{ }}if; \hfill \\
   \h\h\h end{\text{ }}do; \hfill \\
  \h\h\h {\text{coef}}: = {\text{coef}} + w*H_j ; \hfill \\
  \h\h end{\text{ }}do; \hfill \\
  \h\h f_k: = f_k + \kappa (m,n)*\text{coef}*\cos (m\pi x)*\cos (n\pi y); \hfill \\
  \h\h end{\text{ }}do; \hfill \\
  end{\text{ }}do; \hfill \\
\end{gathered}
\]
\\\h In the two following examples, we will show the ill-posedness of the problem and how our regularization treat it.
\\{\bf Example 1.} Consider the exact data
\[
\varphi _0 (t) = \pi ^2 e^{ - 4\pi ^2 t} ,g_0 (x,y) = \left( {1 + \cos (\pi x)} \right)\cos (\pi y),
\]
Then the system (1) has the exact solution
\bqq
  u_0 (x,y,t) &=& e^{ - 4\pi ^2 t} \left( {1 + \cos (\pi x)} \right)\cos (\pi y), \hfill \\
  f_0 (x,y) &=&  - 3\cos (\pi y) - 2\cos (\pi x)\cos (\pi y).
\eqq
For each integer $k\ge 1$, take the disturbed data
\[
\varphi _k (t) = \varphi _0 (t),g_k (x,y) = g_0 (x,y) + \frac{\pi }
{k}\sin ^2 (k\pi x)\cos (k\pi y).
\]
Then the system (1) has the disturbed solution
\bqq
 \widetilde u_k (x,y,t) &=& u_0 (x,y,t) + \frac{\pi }
{k}e^{ - 4\pi ^2 t} \sin ^2 (k\pi x)\cos (k\pi y), \hfill \\
\widetilde f_k (x,y) &=& f_0 (x,y) + \frac{\pi }
{k}\left( {(5k^2  - 4)\sin ^2 (k\pi x) - 2k^2 } \right)\cos (k\pi y).
\eqq
We see that
\[
\left\| {g_k  - g} \right\|_{L^1 (\Omega )}  = \frac{1}
{k},\left\| {\widetilde f_k  - f_0 } \right\|_{L^2 (\Omega )}  = \frac{\pi }
{4}\sqrt {27k^2  - 56 + \frac{{48}}
{{k^2 }}}.
\]
Hence, if $k$ is large then a small error of the data will cause a large error of the solution. Thus the problem is ill-posed and a regularization is necessary.
\\\h Corresponding to the error of the data $\varepsilon=k^{-1}=10^{-2}$, our regularized procedure produces the regularized solution
$$
f_k  (x,y) = - 2.999721\cos (\pi y) - 1.997145\cos (\pi x)\cos (\pi y).
$$
The error between the regularized solution and the exact solution is
$$\left\| {f_k   - f_0 } \right\|_{L^2 (\Omega )}  = 0.001441.$$
In this case, the exact solution has form of a truncated Fourier series, and hence the approximation is very good.
\\{\bf Example 2.} Corresponding to the exact data
\[
\varphi _0 (t) = e^{t} ,~g_0 (x) = \left( {x\cos (1 - x) + \sin (1 - x) - 1} \right)(2y^3  - 3y^2 ),
\]
the exact solution of the system (1) is
\bqq
  u_0 (x,y,t) &=&e^t \left( {x\cos (1 - x) + \sin (1 - x) - 1} \right)(2y^3  - 3y^2 ), \hfill \\
  f_0 (x,y) &=&\left( {2x\cos (1 - x) - 1} \right)(2y^3  - 3y^2 )+\hfill\\
&~&-\left( {x\cos (1 - x) + \sin (1 - x) - 1} \right)(12y - 6).
\eqq
For each integer $k\ge 1$, the disturbed data
\[
\varphi _k (t) = \varphi _0 (t),g_k (x,y) = g_0 (x,y) + \frac{\pi }
{k}\sin ^2 (k\pi x)\cos (2\pi y)
\]
corresponds the disturbed solution
\bqq
 \widetilde u_k (x,y,t) &=&u_0 (x,y,t) + \frac{\pi }{k}e^t \sin ^2 (k\pi x)\cos (2\pi y), \hfill \\
\widetilde f_k (x,y) &=&f_0 (x,y) + \frac{\pi }{k}\left( {(4k^2 \pi ^2  + 4\pi ^2  + 1)\sin ^2 (k\pi x) - 2k^2 \pi ^2 } \right)\cos (2\pi y).
\eqq
We see that
\bqq
\left\| {g_k  - g} \right\|_{L^1 (\Omega )} &=& \frac{1}
{k}\to 0,\hfill\\
\left\| {\widetilde f_k  - f_0 } \right\|_{L^2 (\Omega )}  &=& \frac{\pi }
{4}\sqrt {16\pi ^4 k^2  + 32\pi ^4  + 8\pi ^2  + \frac{{48\pi ^4  + 24\pi ^2  + 3}}
{{k^2 }}} \to +\infty
\eqq
as $k\to +\infty$. Thus a small error of the data causes a large error of the solution.
\text{}
\centerline{
\includegraphics[width=4in]{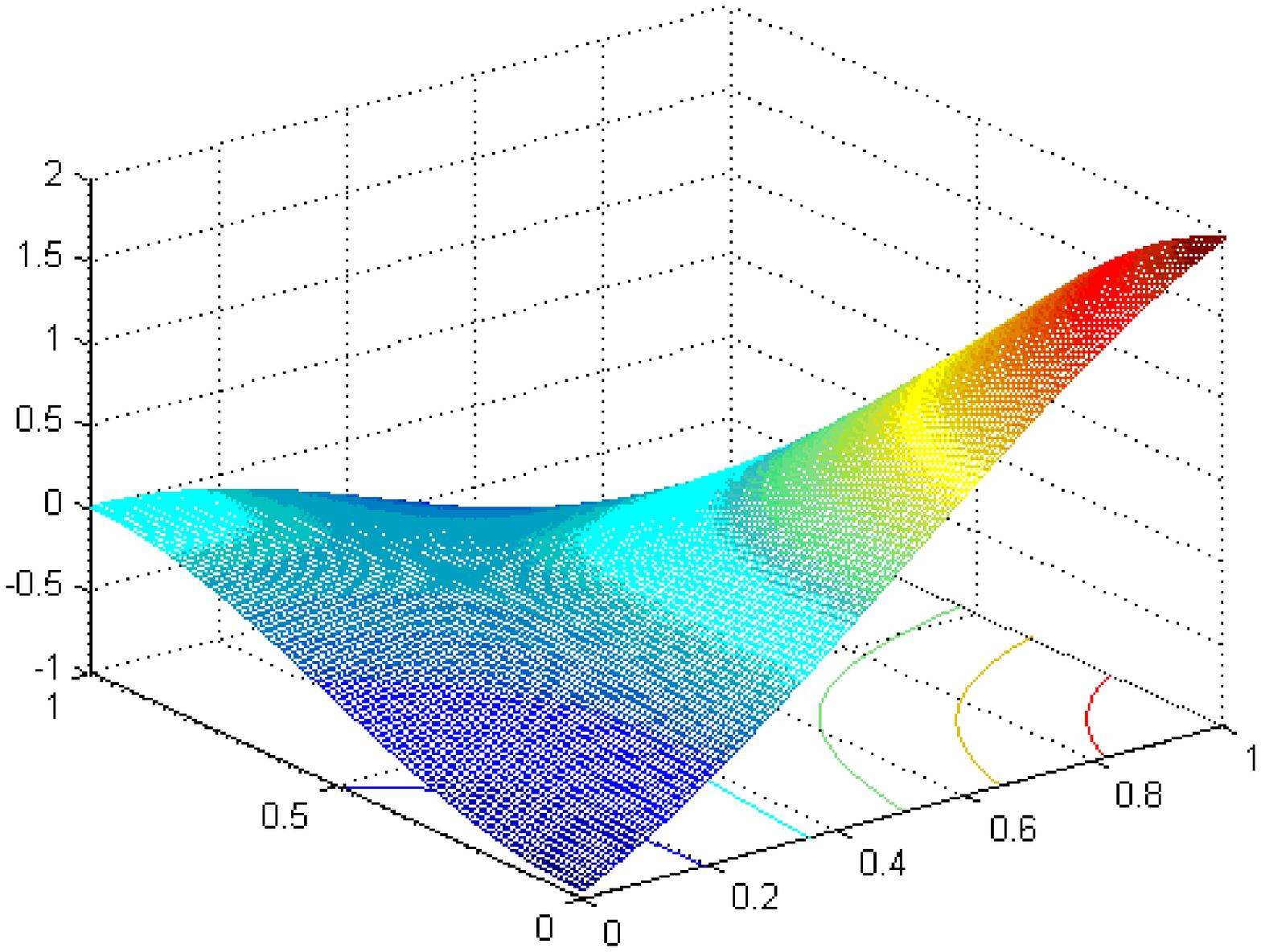}}
\center{Figure 1. The exact solution $f_0(x,y)$.}
\endcenter
\text{}
\centerline{
\includegraphics[width=4in]{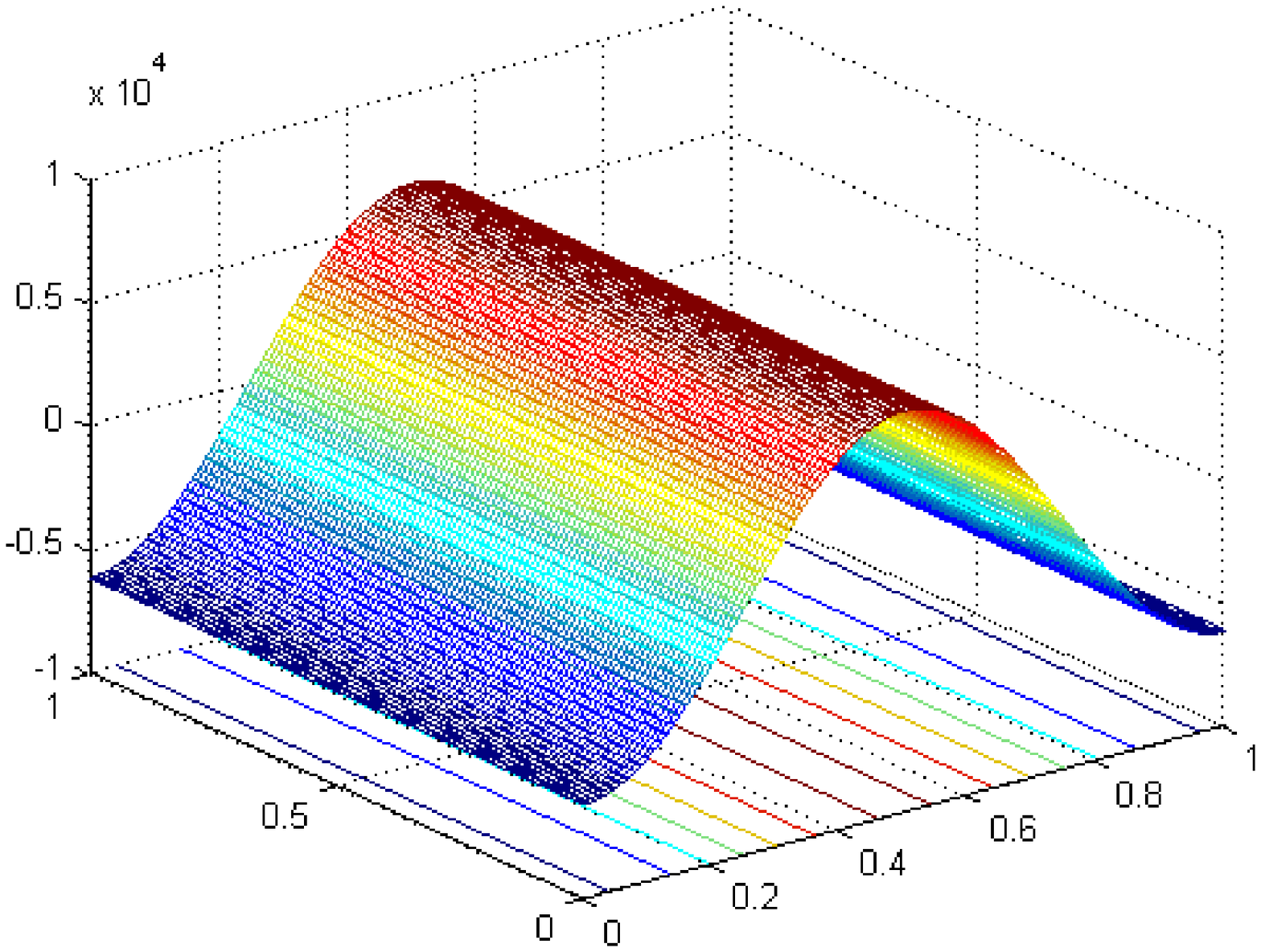}}
\center{Figure 2. The disturbed solution $\widetilde f_k$ with
$k=100$.}\endcenter
 Now, we use our regularized procedure for
$\varepsilon=n^{-1}=10^{-2}$. The regularized solution is \bqq
f_k(x,y)&=&0.040435 + 0.426992\cos (\pi x)\hfill\\
&~&-0.431701\cos (\pi y) - 0.800509\cos (\pi x)\cos (\pi y).
\eqq
\text{}
\centerline{
\includegraphics[width=4in]{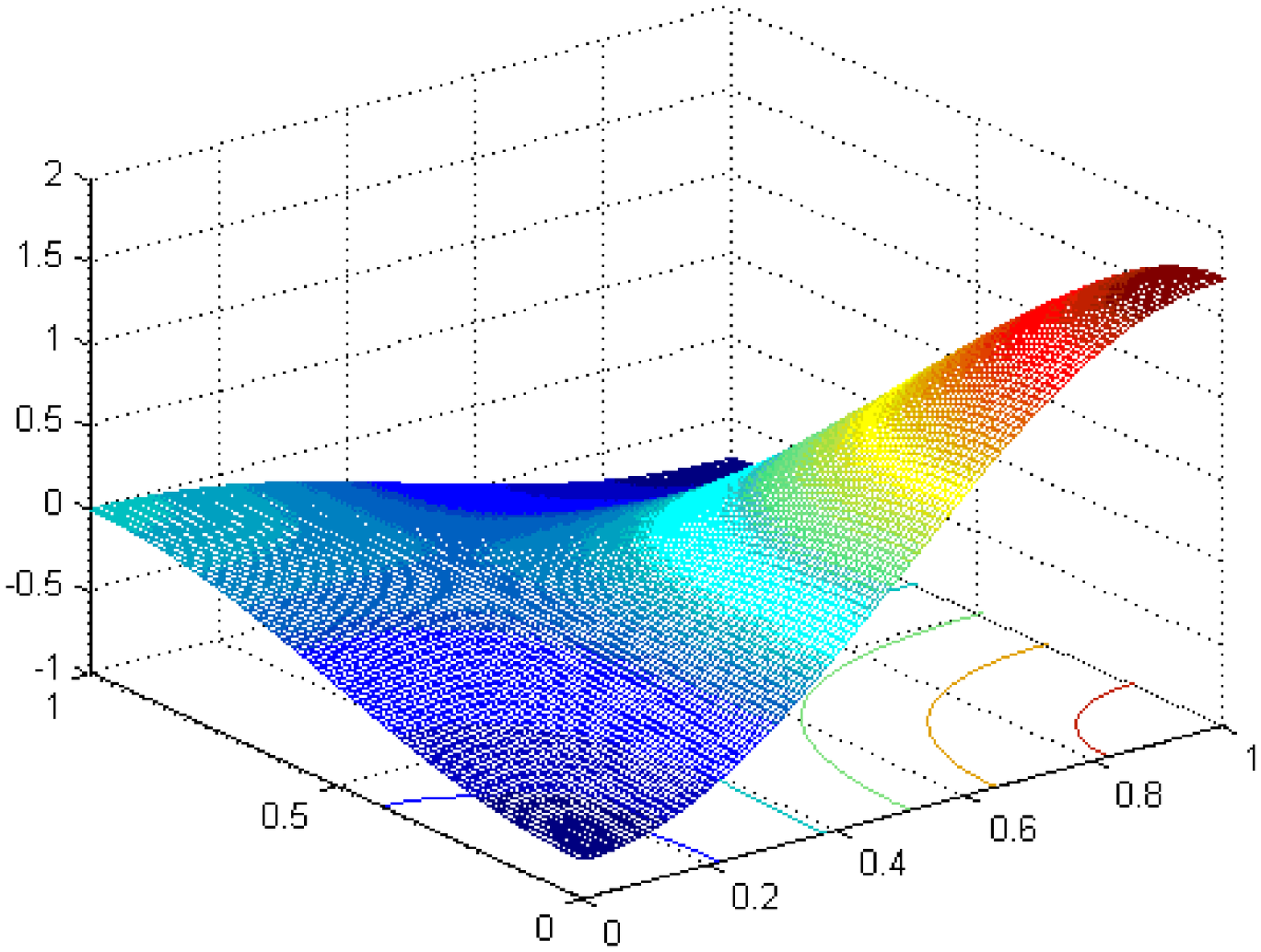}}
\center{Figure 3. The regularized solution $f_k$ with $k=100$.}
\endcenter
\h The error between the regularized solution and the exact solution is $\left\| {f_{k}  - f_0 } \right\|_{L^2 (\Omega)}= 0.059997$. To see the effect of the regularization, we note that the disturbed solution corresponding to $k=100$ causes a so large error $\left\| {\widetilde f_k  - f_0 } \right\|_{L^2 (\Omega )}=1.24\times 10^6$. Figure 1, Figure 2 and Figure 3 give a visual comparison between the exact solution, the disturbed solution and the regularized solution in the second example.


\begin{thebibliography}{17}
\bibitem{C68} J.R.Cannon, {\it Determination of an unknown heat source from overspecified boundary data}, SIAM J. Numer. Anal. 5(1968) 275-286.
\bibitem{CP90} J.R.Cannon, S.Peùrez Esteva, {\it Some stability estimates for a heat source in terms of over specified data in the 3-D heat equation}, J. Math. Anal. Appl. 147(1990), no.2, 363-371.
\bibitem{CP91} J.R.Cannon, S.Peùrez Esteva, {\it Uniqueness and stability of 3D heat source}, Inverse Problems 7(1991), no.1, 57-62.
\bibitem{CY04} M. Choulli, M. Yamamoto, {\it Conditional stability in determining a heat
source}, J. Inverse Ill-Posed Problems 12(2004), no.3, 233-243.
\bibitem{FL06} A.Farcas, D.Lesnic, {\it The boundary-element method for the determination of a heat source dependent on one variable}, J. Engrg. Math. 54(2006), no.4, 375-388.
\bibitem{Ya} B.Ya. Levin, {\it Lectures on Entire Functions}, Trans Math Monographs, Vol.150, AMS, Providence, Rhole Island, 1996.
\bibitem{STY02} S.Saitoh, V.K.Tuan, M.Yamamoto, {\it Reverse convolution inequalities and applications to inverse heat source problems}, JIPAM. J. Inequal. Pure. Appl. Math. 3(2002), no.5, Article 80 (electronic).
\bibitem{STY03} S.Saitoh, V.K.Tuan, M.Yamamoto, {\it Convolution inequalities and applications}, JIPAM. J. Inequal. Pure. Appl. Math. 4(2003), no.3, Article 50 (electronic).
\bibitem{SZN05} A.Shidfar, A.Zakeri, A.Neisi, {\it A two-dimensional inverse heat conduction problem for estimating heat source}, Int. J. Math. Math. Sci. (2005), no.10, 1933-1941.
\bibitem{TLD05} D.D.Trong, N.T.Long, P.N.Dinh Alain, {\it Nonhomogeneous heat equation: Identification and regularization for the inhomogeneous term}, J. Math. Anal. Appl. 312 (2005), 93-104.
\bibitem{TQD06} D.D.Trong, P.H.Quan, P.N.Dinh Alain, {\it Determination of a two-dimensional heat source: Uniqueness, regularization and error estimate}, J. Comput. Appl. Math. 191(2006), no.1, 50-67.
\bibitem{WZ99} P.Wang, K.Zheng, {\it Determination of the source/sink term in a heat equation}, Fourth Mississippi State Conference on Differential Equations and Computational Simulations, Electronic J. Diff. Equations, Conference 03(1999) 119-125.
\bibitem{WZ06} P.Wang, K.Zheng, {\it Reconstruction of spacial heat sources in heat conduction problems}, Appl. Anal. 85(2006), no.5, 459-465.
\bibitem{Y93} M. Yamamoto, {\it Conditional stability in determination of force terms of heat
equations in a rectangle}, Math. Comput. Modelling 18 (1993), no.1, 79-88.
\bibitem{Y94} M. Yamamoto, {\it Conditional stability in determination of densities of heat source in a bounded domain}, Control and estimation of distributed parameter systems: nonlinear phenomena (Vorau, 1993), 359-370, Internat. Ser. Numer. Math., 118, Birkhauser, Babel, 1994.
\bibitem{YZ01} M. Yamamoto, J.Zou, {\it Simultaneous reconstruction of the initial temperature and heat radiative coefficient}, Inverse Problems 17(2001), no.4, 1181-1202.
\end{thebibliography}
\end{document}